\documentclass[12pt]{report}
\usepackage[cp1251]{inputenc}
\usepackage[english,russian,ukrainian]{babel}
\usepackage{latexsym,amsfonts,amsthm,amssymb,amsmath}
\usepackage{graphicx,graphics,hhline}
\usepackage{euscript}
\usepackage[all]{xy}
\begin{document}

\Large

\noindent УДК 513.83; 517.5
\vskip 2mm
\smallskip

\noindent{\bf T.M.~Osipchuk (Т.М.~Осіпчук)},  otm82@mail.ru
\vskip 2mm

\noindent Institute of Mathematics NAS of Ukraine

\vskip 4mm

\noindent{\bf ON THE SHADOW PROBLEM FOR DOMAINS  

\noindent   IN THE EUCLIDEAN SPACES }
\vskip 5mm
\noindent{\bf ЗАДАЧА ПРО ТІНЬ ДЛЯ ОБЛАСТЕЙ

\noindent  В  ЕВКЛІДОВИХ ПРОСТОРАХ}

\vskip 8mm

\small

\noindent{In the present work, the problem about shadow, generalized on domains of space $\mathbb{R}^n$, $n\le 3$, is investigated. Here the shadow problem means to find the minimal number of balls satisfying some conditions an such that every line passing through the given point intersects at least one ball of the collection. It is proved that to generate the shadow at every given point of any domain of the space $\mathbb{R}^3$ ($\mathbb{R}^2$) with collection of  mutually non-overlapping closed or open balls which do not hold the point and with centers on  the boundary of the domain, it is sufficient to have four balls.}

\vskip 3mm

\noindent{В роботі досліджується задача про тінь, узагальнена на області простору $\mathbb{R}^n$, $n\le 3$. Під задачею про тінь тут мається на увазі знаходження мінімального числа куль, що задовільняють деяким умовам, і таких, що кожна пряма, що проходить через дану точку, перетне принаймні одну кулю з набору.  Доведено, що для того, щоб створити тінь в кожній даній точці довільної області простору $\mathbb{R}^3$ ($\mathbb{R}^2$) набором замкнених або відкрити куль, які попарно не перетинаються, не містять дану точку та з центрами на межі області,  достаточно чотирьох (двох) таких куль. }

\vskip 3mm

\noindent{В работе исследуется задача о тени, обобщенная на области пространства $\mathbb{R}^n$, $n\le 3$. Под задачей о тени подразумевается нахождение минимального количества шаров, удовлетворяющих некоторым условиям, и таких, что каждая прямая, проходящая через заданную точку, пересечет хотя бы один шар из набора. Доказано, что для того, чтобы создать тень в каждой заданной точке произвольной области пространства $\mathbb{R}^3$ ($\mathbb{R}^2$) набором замкнутых или открытых шаров, попарно не пересекающихся, не содержащих заданную точку и с центрами на границе области, достаточно четырех (двух) таких шаров.}

\newpage

\large

\begin{center} {\bf Introduction} \end{center}

In 1982 G. Khudaiberganov proposed the problem about shadow \cite{Hud} that can be formulated as follows:  {\it what is the minimal number of closed (open) non-overlapping balls in the $n$-dimensional Euclidean space with radiuses less than the sphere radius,  with centers on it and such that every line passing through the center of the sphere would intersect at least one of the balls.}

This problem (and similar to it) may be written  briefly as: {\it what minimal number of some objects generates the shadow at the point.}

This problem was solved by G. Khudaiberganov for the case $n=2$: it was proved that two balls are sufficient for a circumference on the plane. For all that, it was also made the assumption that for the case $n>2$ the minimal number of such balls is exactly equal to $n$.

Yu. Zelinskii mentioned this problem in many of his works (cf., e.g.  [1\,-\,2]). The problem is also interesting from the standpoint of convex analysis as a special case of the problem about membership of a point to the 1-hull of the union of some collection of balls. In \cite{Zel3}, Yu. Zelinskii and his students proved that three balls are not sufficient for the case $n=3$, but it is possible to generate the shadow at the center of the sphere with four balls. In their work  it is also proved that for the general case the minimal number is $n+1$ balls.   Thus, G. Khudaiberganov's assumption was wrong.  In \cite{Zel3}, it is also proposed another method of solving the problem for the case $n=2$ which gives some numerical estimates.

In \cite{Zel3_1} it is proved that in order that collection of mutually non-overlapping closed (open) balls of $\mathbb{R}^n$, $n\ge 2$ which do not hold the chosen point of  $\mathbb{R}^n$ generate the shadow at the point, it is necessary and sufficient to have $n$ balls.

All results mentioned above and some others related to shadow problem can also be found in the survey article \cite{Zel3_2}.

The following question  is naturally arising. What can be done in this problem if we will consider other surfaces instead of sphere? In  \cite{Osipchuk_Tkach2}  (cf., e.g.  \cite{Osipchuk_Tkach1}) there were considered prolate spheroids with value of ratio between its major and minor axis more then $ 2 \sqrt2$. Authors proved that in order that the union of mutually non-overlapping closed (open)  balls which do not hold the ellipsoid center and with centers on  the ellipsoid generate the shadow at the ellipsoid center, it is necessary and sufficient to have three balls. But the question still was opened for the other prolate ellipsoids. In the present work the shadow problem for every domain and every its point is solved,  which also gives the answer for the ellipsoids.

\begin{center} {\bf Statement of Problem } \end{center}

\noindent{\bf Lemma 1.}  {\it Let two open or closed non-overlapping balls of $\mathbb{R}^n$,  with centers on a sphere $S^{n-1}$ and with radiuses less than the sphere radius, are given. Then, every ball with center outside of the sphere, homothetic to the smaller ball (homothetic to one of the balls, if they are equal) relative to the sphere center, would not intersect the other ball.}

\vskip 2mm

\noindent{\bf Proof.} Without loss of generality, let us consider the unit sphere. We denote the balls as $B_1(O_1,R)$, $B_2(O_2,r)$ with radiuses $R$, $r$ and centers $O_1$, $O_2$ respectively and let \begin{equation}\label{1}r\le R.\end{equation} Let the distance between the balls' centers is equal to $d$, then, the condition such that the balls are non-overlapping is equivalent to the following  inequality
\begin{equation}\label{2}
r+R\le d.
\end{equation}

Further, we will consider the angle $\alpha$ (figure 1) between the segment $OO_2$ joining the sphere center and the center of the ball $B_2$ and the ray coming out of the sphere center and which is tangent to the ball $B_2$ at the point $A$. And we will also consider the angle $\varphi$ between the segment $OO_2$ and the perpendicular to $d$. We do not consider the case when the triangle $O_1OO_2$ is degenerated into the segment $O_1O_2$, as it coincides to diametrically opposite situation of the balls $B_1$, $B_2$ which obeys the conditions of lemma 1. It is easy to see that $\angle\alpha\le\angle\varphi$.  Indeed, let us consider right triangles $OAO_2$ and $ODO_2$. We get from them $\sin\alpha=r$, $\sin\varphi=\dfrac{1}{2}d$. The inequalities  (\ref{1}), (\ref{2}) give $r\le \dfrac{1}{2}d$, that is why $$\sin\alpha\le\sin\varphi,$$ which leads to $\angle\alpha\le\angle\varphi$.

\begin{center}
\includegraphics[width=10 cm]{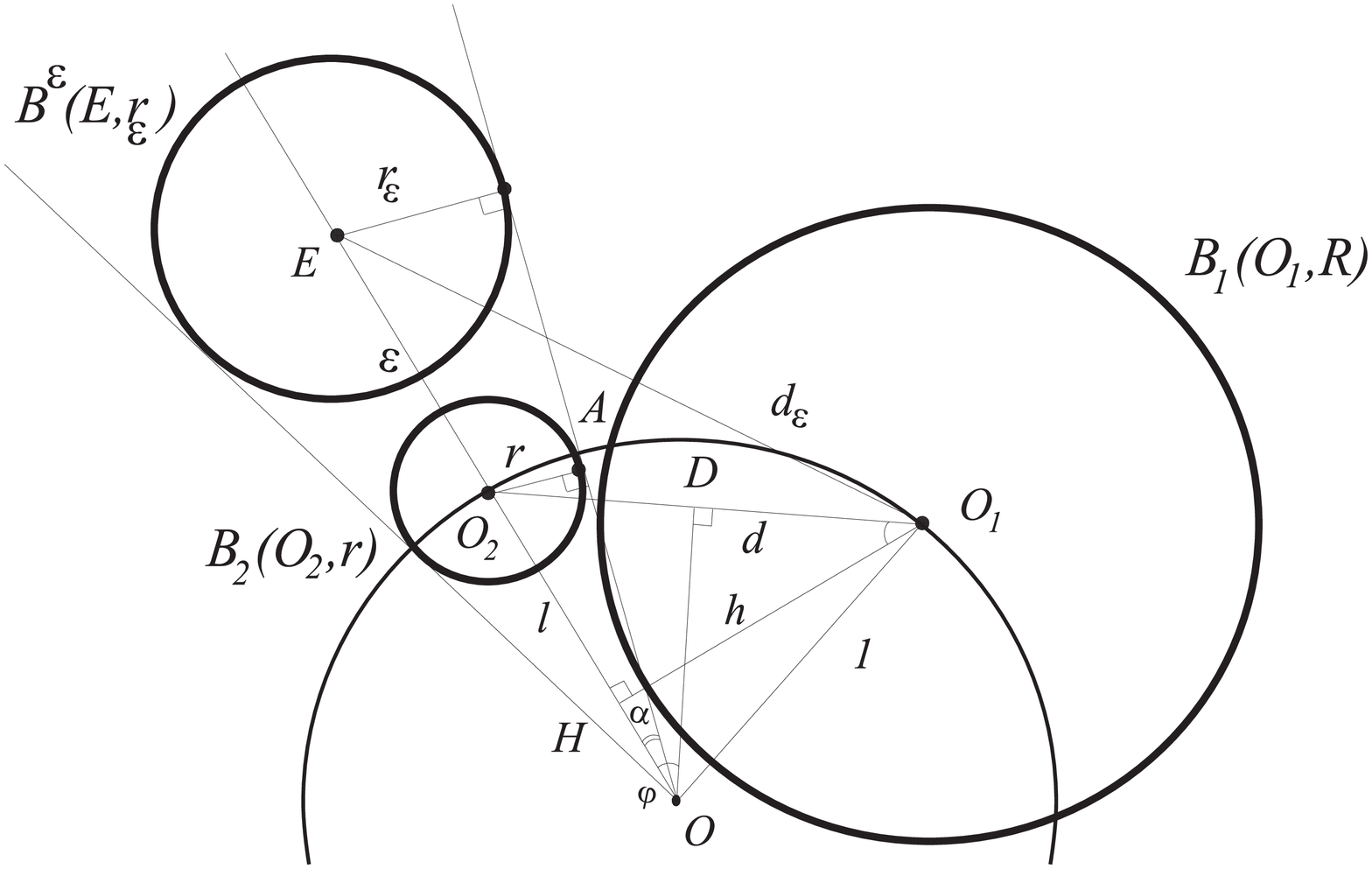}
\end{center}
\small\begin{center} Figure 1 \end{center}
\large

Now let us consider a ball $B^\varepsilon(E,r_\varepsilon)$ with the center outside of the sphere at the point $E$, situated at the distance $\varepsilon$ from the center of the ball $B_2$, with radius $r_\varepsilon$ and which is homothetic to the ball $B_2$ relative to the sphere center $O$. Let us show that for this ball  it is also valid the inequality $$
r_\varepsilon+R\le d_\varepsilon,
$$
where  $d_\varepsilon$ be the distance between the centers of the balls $B^\varepsilon$, $B_1$.

First we notice that  the equality $r_\varepsilon=r+\varepsilon\sin\alpha$ and  (\ref{2}) give
$$
r_\varepsilon+R\le d+\varepsilon\sin\alpha.
$$
Further, we drop the altitude $O_1H=h$ in the triangle  $\triangle O_1OO_2$ from the vertex  $O_1$, then $\angle O_2O_1O=\varphi$ and $$\label{3}HO_2=l=d\sin\varphi.$$
From $\triangle O_1HO_2$ we get
$$
d^2=h^2+l^2,
$$
and from $\triangle O_1HE$ and two privies formulas we have
\begin{equation}\label{5}
d^2_\varepsilon=h^2+(l+\varepsilon)^2=d^2+\varepsilon^2+2l\varepsilon=d^2+\varepsilon^2+2d\varepsilon\sin\varphi.
\end{equation}
Let us consider the formula
Розглянемо вираз
\begin{equation}\label{6}
(d+\varepsilon\sin\alpha)^2=d^2+\varepsilon^2\sin^2\alpha+2d\varepsilon\sin\alpha.
\end{equation}
Now from (\ref{5}),  (\ref{6}) it is easy to see that
$$
d+\varepsilon\sin\alpha\le d_\varepsilon,
$$
so, $r_\varepsilon+R\le d_\varepsilon$, which was to be demonstrated. Lemma 1 is proved.

\vskip 2mm

It is not difficult to give a counter example with balls, homothetic to the smaller ball relative to the sphere center and with centers inside the sphere.

\vskip 2mm

\noindent{\bf Corollary 1.}  {\it Let two open or closed non-overlapping balls of $\mathbb{R}^n$,  with centers on a sphere $S^{n-1}$ and with radiuses less than the sphere radius, are given. Then, every ball with center inside of the sphere, homothetic to the bigger ball (homothetic to one of the balls, if they are equal) relative to the sphere center, would not intersect the other ball.}

\vskip 2mm

\noindent{\bf Corollary 2.}  {\it Let two open or closed non-overlapping balls of $\mathbb{R}^n$,  with centers on a sphere $S^{n-1}$ and with radiuses less than the sphere radius, are given. Then, every ball, homothetic to the smaller ball (homothetic to one of the balls, if they are equal), relative to the sphere center, with coefficient of homothety $k_1$ would not intersect every ball, homothetic to the bigger one, relative to the sphere center, with coefficient of homothety $k_2$, if $k_1\ge k_2$.}

 \vskip 2mm
It was proved in \cite{Zel3} that it is possible to generate the shadow in the center of sphere $S^{n-1}$ by four given balls of  $\mathbb{R}^n$, proposing one of the placement methods of their centers on the sphere. In the following lemma we give our method, convenient for further statements. We will also  consider double-napped cones with common vertex at the center of the sphere with generators tangent to the balls. Here and further such cones will be referred to as {\it cones under the balls}.

  \vskip 2mm

\noindent{\bf Lemma 2.}  {\it  There are four closed (opened) non-overlapping balls of $\mathbb{R}^n$,  with centers in the vertexes of regular triangular pyramid and with radiuses less than the radius of  sphere circumscribed about  the pyramid such that generate the shadow in the sphere center.}

\vskip 2mm

\noindent{\bf Proof.} Without loss of generality, we will consider the unit sphere $S^{2}$ which center for convenience will be placed at the point  $O_1=(0,-1,0)$ (Figure 2). Let us place the center of the unit open ball $B_1(O,1)$ at the origin,  $O=(0,0,0)$. Then the ball $B_1$ will generate shadow for the sphere center everywhere except sphere equatorial plane $\Sigma=\{(x,y,z): y=-1\}$ tangent to the ball $B_1$ at the point $O_1$. Let us establish dependence of the angle $\varphi$ at which the intersection of the plane $\Sigma$ with every ball $B^X(X,r(X))$ tangent to the ball $B_1$ with center at point $X\in\mathbb{R}^3$ and radius $r$ can be seen from the point $O_1$, on the point $X$.

\begin{center}
\includegraphics[width=10 cm]{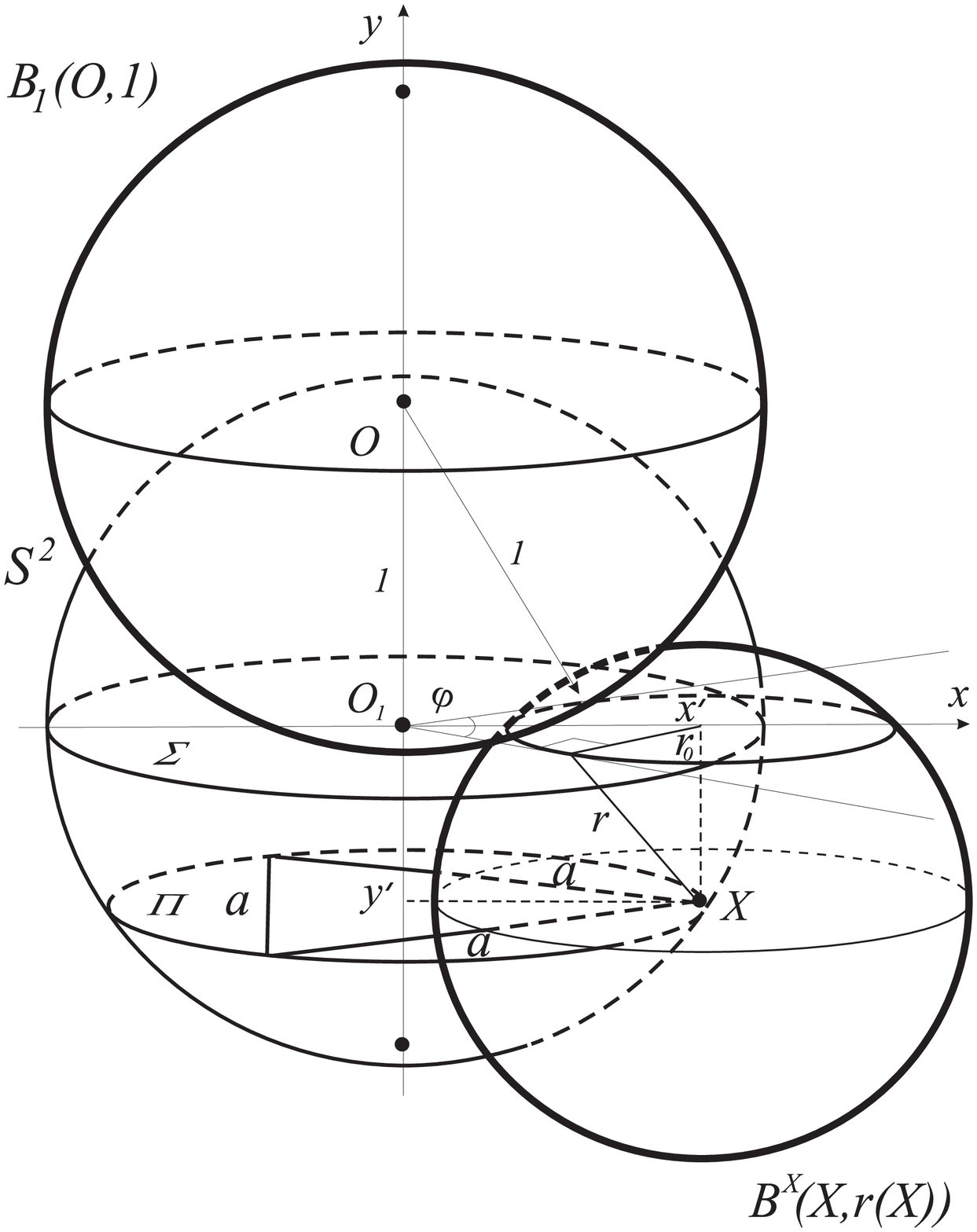}
\end{center}
\small\begin{center} Figure 2 \end{center}
\large

Without loss of generality, let us consider opened balls  $B^X$ with centers at points $X=(x,y,0)$ of the half-plane  $xOy$, $x> 0$, and such that have nonempty intersection with plane $\Sigma$. Then this intersection will be an opened ball with center at a point $(x,-1,0)$ and some radius $r_0$. We obtain:
$$
\sin\frac{\varphi}{2}=\frac{r_0}{x},
$$
$$
r_0=\sqrt{r^2-(-y-1)^2},
$$
$$
r=\sqrt{x^2+y^2}-1.
$$
Dropping intermediate calculations, we get:
$$
\sin\frac{\varphi}{2}=\frac{\sqrt{x^2-2y-2\sqrt{x^2+y^2}}}{x},\,\,\, x> 0.
$$
\begin{equation}\label{8}
\frac{1}{4}x^2cos^2\frac{\varphi}{2}-\frac{1}{cos^2\frac{\varphi}{2}}=y.
\end{equation}

Further, we will consider among all balls $B^X$ those, who's centers are on the sphere. The following system is obtained:
$$
\left\{\begin{array}{l}
\sin\dfrac{\varphi}{2}=\dfrac{\sqrt{x^2-2y-2\sqrt{x^2+y^2}}}{x},\\
x^2+(y+1)^2=1,\\
x> 0,
\end{array}\right.
$$
that is equivalent to equation
\begin{equation}\label{7}
\sin^2\frac{\varphi}{2}=\frac{-y^2-4y-2\sqrt{-2y}}{-y^2-2y}, \,\,\, \mbox{при}\,\,\, -2 < y<0,
\end{equation}
or to equation
$$
\cos^2\dfrac{\varphi}{2}y^2+2\left(\cos^2\dfrac{\varphi}{2}+2\right)y+\dfrac{4}{\cos^2\dfrac{\varphi}{2}}=0 \,\,\, \mbox{при}\,\,\, -2 < y<0.
$$
Let us find the value of the coordinate $y$ of the center of the ball $B^X$ that overlaps $\angle\varphi=\dfrac{\pi}{3}$. We obtain:

$$
y'=\frac{1}{3}(\sqrt{57}-11)\approx-1,15.
$$
At this point
$$
x'=\frac{4}{3}\sqrt{\sqrt{57}-7}\approx0,99.
$$
Thus, every ball, tangent to the ball $B_1$ with center on the sphere at the level $y=y'$, will overlap a sector with $\angle\varphi=\dfrac{\pi}{3}$ in the plane $\Sigma$. Herewith, it is easy to see that radius of such ball $r'=\sqrt{(x')^2+(y')^2}-1\approx 0,52$ is less than a half of side $a=\sqrt{3}x'\approx 1,71$ of regular triangle inscribed in circumference on the sphere in the plane $\Pi=\{(x,y,z): y=y'\}$. Thus, three balls, tangent to the ball $B_1$ with centers at the vertexes of the triangle, will not overlap.

Now, let us consider plane $\Pi_\varepsilon=\{(x,y,z): y=y'-\varepsilon\}$. We place the centers of three balls tangent to the ball $B_1$ at the vertexes of a regular triangle inscribed in circumference that is the intersection of plane $\Pi_\varepsilon$ with sphere. At this point, we choose $\varepsilon>0$ small enough for these balls to not intersect yet. Since the function $\sin^2\dfrac{\varphi}{2}$ in (\ref{7}) is decreasing, then every such a ball will overlap sector $\angle\varphi>\dfrac{\pi}{3}$ in the plane $\Sigma$. Thus, altogether they will generate the shadow at the sphere center but only in the plane $\Sigma$ and in the surface of every cone with vertex at the point $O_1$ and with apex angle close enough to $\pi$. Now, we will decrease the radius of the ball $B_1$  in order that the angle of the cone under it remain close enough to $\pi$. And, thus, we obtain set of balls satisfying conditions of the lemma 2.

The lemma 2 is proved.

 \vskip 2mm

 In \cite{Zel2}, there are considered sets similar to the convex ones which are defined as follows.

\noindent {\bf Definition 1.} (\cite{Zel2}) {\it The set  $E\subset\mathbb{R}^n$  is {\it $m$-convex relative to the point }  $x\in \mathbb{R}^n\setminus E$, if there exists an  $m$-dimensional plane $L$ such that $x\in L$  and  $L\cap E=\emptyset$.}

\noindent {\bf Definition 2.} (\cite{Zel2}) {\it The set $E\subset\mathbb{R}^n$  is {\it $m$-convex} if it is $m$-convex relative to every point $x\in \mathbb{R}^n\setminus E$.}

It is easy to verify that both definitions satisfy the axiom of convexity: the intersection of every subfamily of such sets also satisfies the definition. For every set $E\subset\mathbb{R}^n$, it can be considered the minimal $m$-convex set containing $E$  which is called {\it $m$-hull} of the set $E$.

Now, the shadow problem can be considered as a partial case of the membership of a point to the 1-hull of a union of some collection of balls.

Further, let domain be an open connected set.

\vskip 2mm

\noindent {\bf Theorem 1.}  {\it In order that every given point $x_0$ of a domain $D\subset\mathbb{R}^2$ belong to $1$-hull of  mutually non-overlapping closed (open) balls which do not hold the point $x_0$ and with centers on  the boundary of the domain $D$, it is sufficient to have two balls.}

 \vskip 2mm

 \noindent {\bf Proof.} Let us consider the circumference $S^1$ with center at the point $x_0$ and with maximal radius such that its interior is contained in domain  $D$. Then, let us construct the collection of two mutually non-overlapping closed (open) circles which do not hold the point $x_0$ and with centers on the circumference as it was done in \cite{Hud}. Along with, let us situate the center of the bigger circle at the point $x\in\partial D\cap S^1$. Now, we construct the circle that is homothetic to the smaller one relative to the point  $x_0$ and with center on $\partial D$. Then, due to lemma 1, this ball do not intersect the bigger one. The theorem is proved.

\vskip 2mm

\noindent {\bf Theorem 2.}  {\it In order that every given point $x_0$ of a domain $D\subset\mathbb{R}^3$ belong to $1$-hull of  mutually non-overlapping closed (open) balls which do not hold the point $x_0$ and with centers on  the boundary of the domain $D$, it is sufficient to have four balls.}

 \vskip 2mm

 \noindent {\bf Proof.} Let us consider sphere $S^2$ with center at the point $x_0$ and with maximal radius such
 that its interior is contained in the domain  $D$. Then, let us construct the collection of four mutually
 non-overlapping closed (open) balls which do not hold the point $x_0$ and with centers on the sphere as it was
 done in lemma 2. Along with, let us situate the center of the biggest ball at the point $x\in\partial D\cap S^2$.
 Now, we use the homothety to the rest three balls relative to the point  $x_0$  in such a way as the centers of
 homothetic balls to be on $\partial D$.
 % in such a way as to (do something)
%(A disguised message is included in such a way as to be audible only when the disc is spun backwards.)
%so that, so as
Then, due to lemma 1 each of these balls would not intersect the biggest one and due to corollary  2 these balls would not intersect each other. The theorem is proved.

 \vskip 2mm

 There exist  domain of  $\mathbb{R}^3$ with  point that belongs to $1$-hull of three  mutually non-overlapping closed (open) balls which do not hold the point and with centers on  the boundary of the domain. A prolate ellipsoid with value of ratio between its major and minor axis more then $ 2 \sqrt2$ is the example of such domain and the ellipsoid center is the point (\cite{Osipchuk_Tkach1}, \cite{Osipchuk_Tkach2}).

\renewcommand{\bibname}{References}

\end{document}